\documentclass[11pt]{article}
\usepackage[utf8]{inputenc}
\usepackage[a4paper, margin = 0.5in, bottom = 1in]{geometry}
\usepackage{soul}

\usepackage[utf8]{inputenc}
\usepackage[a4paper, margin = 0.5in, bottom = 1in]{geometry}
\usepackage{soul}
\usepackage{authblk}
\usepackage{minted}
\usepackage{relsize}

\date{}
\usepackage{amsmath}
\usepackage{amsfonts}
\usepackage{amssymb}
\usepackage{mathtools}
\usepackage{mathrsfs}
\usepackage{bm}
\usepackage{wasysym}
\usepackage{graphicx}
\usepackage{caption}
\usepackage{cite}
\usepackage{subcaption}
\usepackage{skak}
\usepackage[table]{xcolor}

\usepackage{amsfonts}
\usepackage{subcaption}
\usepackage{lmodern}
\usepackage{xcolor}
\usepackage{csquotes}
\usepackage{mathtools}
\usepackage{float}
\usepackage{amsmath,amssymb,amsthm,comment}

\usepackage{authblk}
\usepackage{minted}

\date{}
\usepackage{amsmath}
\usepackage{amsfonts}
\usepackage{amssymb}
\usepackage{mathtools}
\usepackage{mathrsfs}
\usepackage{bm}
\usepackage{wasysym}
\usepackage{graphicx}
\usepackage{caption}
\usepackage{cite}
\usepackage{subcaption}
\usepackage{skak}
\usepackage[table]{xcolor}

\usepackage{amsfonts}
\usepackage{subcaption}
\usepackage{lmodern}
\usepackage{xcolor}
\usepackage{csquotes}
\usepackage{mathtools}
\usepackage{float}
\usepackage{amsmath,amssymb,amsthm,comment}
\usepackage{amsfonts}
\usepackage{subcaption}
\usepackage{xcolor}
\usepackage{csquotes}
\usepackage{mathtools}
\usepackage{float}
\usepackage{amsmath,amssymb,amsthm}
\usepackage[colorlinks,citecolor=red,urlcolor=blue,linkcolor=blue]{hyperref}
\usepackage{mathrsfs,dsfont}
\usepackage{graphicx}
\usepackage{enumerate}
\usepackage{nicefrac}
\usepackage[font={normalsize},width=1.12\linewidth]{caption}
\usepackage{comment}
\usepackage{hyperref}
\usepackage{pgfplots}
\pgfplotsset{compat=1.15}
\usepackage{mathrsfs}
\usetikzlibrary{arrows}

\newtheorem{nummer}{ }[section]
\newtheorem{thm}[nummer]{Theorem}

\newtheorem{lem}[nummer]{Lemma}
\newtheorem{cor}[nummer]{Corollary}

\newtheorem{defi}[nummer]{Definition}
\newtheorem{conj}[nummer]{Conjecture}

\makeatletter
\def\opargproof[#1]{\par\noindent {\bf #1 }}
 \makeatother

\begin{document}
\medskip\medskip
\begin{center}
\vspace*{50pt}
{\LARGE\bf The Greedy Algorithm for Dissociated Sets}

\bigskip
{\small Sayan Dutta}\\[1.2ex] 
 {\scriptsize 
sayandutta345@gmail.com\\
\href{https://sites.google.com/view/sayan-dutta-homepage}{https://sites.google.com/view/sayan-dutta-homepage}}
\\[3ex]

\end{center}

\begin{abstract}
A set $\mathcal S\subset \mathbb N$ is said to be a \textit{subset-sum-distinct} or \textit{dissociated} if all of its finite subsets have different sums. Alternately, an equivalent classification is if any equality of the form
$$\sum_{s\in \mathcal S} \varepsilon_s \cdot s =0$$
where $\varepsilon_s \in \{-1,0,+1\}$ implies that all the $\varepsilon_s$'s are $0$. For a dissociated set $\mathcal S$, we prove that for $c_\ast = \frac 12 \log_2 \left(\frac \pi 2\right)$ and any $c_\ast-1<C<c_\ast$, we have
$$\mathcal S(n) \,:=\, \mathcal S\cap [1,n] \,\le\, \log_2 n +\frac 12 \log_2\log_2 n + C$$
for all $n\in \mathcal N_C$ with asymptotic density $\mathbf d\left(\mathcal N_C\right)=2-2^{c_\ast-C}$. Further, we consider the greedy algorithm for generating these sets and prove that this algorithm always eventually doubles. Finally, we also consider some generalizations of dissociated sets and prove similar results about them.
\end{abstract}

\section{Introduction}
An $m$-element set $\mathcal S \subset \mathbb N$ is said to have distinct subset sums if all the $2^m$ subsets of $\mathcal S$ have different sums \cite{er1, er2, er3, er4, er5, guy2, stefan, tom, conway, lunnon, Maltby}. These sets have often been called \textit{dissociated sets} in literature where they have been studied in more general settings \cite{diss1, diss2, diss3, diss4, diss5, diss6, diss7}. Although this note is only concerned with subsets of $\mathbb N$, this is the terminology that will be used. An infinite subset of $\mathbb N$ is called dissociated if all its finite subsets are. Another way to describe these sets is
\begin{defi}
    A set $\mathcal S\subset \mathbb N$ is called \textit{dissociated} if any equality of the form
    $$\sum_{s\in \mathcal S} \varepsilon_s \cdot s =0$$
    where $\varepsilon_s \in \{-1,0,+1\}$ implies that all the $\varepsilon_s$'s are $0$.
\end{defi}

An example of such sets is the set of powers of $2$.

Now, for such a set $\mathcal S$, let us define $\mathcal S(n):=\mathcal S \cap [1,n]$ to be the elements of $\mathcal S$ that are less than or equal to $n$. The first observation that was made in this field is the following.

\begin{thm}[Erd\H{o}s]\label{thm:dens}
    For a dissociated set $\mathcal S$, we have $\left\vert\mathcal S(n)\right\vert \le \log_2 n + \mathcal{O}( \log_2\log_2 n)$.
\end{thm}
\begin{proof}
    Consider the map
    $$f : \mathcal P\left(\mathcal S(n)\right) \to \mathbb N$$
    that takes a subset of $\mathcal S(n)$ to the sum of its elements. Since $\mathcal S$ is dissociated, this map is injective. So, $2^{\left\vert\mathcal S(n)\right\vert} \leq n\cdot \left\vert\mathcal S(n)\right\vert +1$, hence completing the proof.
\end{proof}

Theorem \ref{thm:dens} was later improved to
$$\left\vert\mathcal S(n)\right\vert \,<\, \log_2 n \,+\, \frac 12 \log_2\log_2 n \,-\, \log_2 c^\ast \,+\, o(1)$$
where the constant $c^\ast$ has been improved subsequently to $1/4$ (Erd\H{o}s and Moser \cite{er1}), $2\cdot(27)^{-1/2}$ (Alon and Spencer \cite{alone}), $\pi^{-1/2}$ (Elkies \cite{elk}), $3^{-1/2}$ (Bae \cite{bae}, and Guy \cite{guy}), $\sqrt 3\cdot (2\pi)^{-1/2}$ (Aliev \cite{Ali}) and finally $\sqrt 2\cdot \pi^{-1/2}$ (Dubroff, Fox, and Xu \cite{xu}).

A classic conjecture of Erd\H{o}s states
\begin{conj}[Erd\H{o}s, \$ 500]
    For a dissociated set $\mathcal S$, we have
    $$\left\vert\mathcal S(n)\right\vert \leq \log_2 n + c$$
    where $c$ is a fixed constant.
\end{conj}
This problem is listed as Problem $\# 1$ in Bloom's database of Erd\H{o}s problems \cite{conj}.

This note is mainly concerned with the greedy algorithm for generating these sets. It seems that this question hasn't been discussed in literature before, although there are some interesting observations to be made. But before that, we will give a slightly generalized version of the theorem in \cite{xu} using a result of Harper \cite{har}. The idea is to show that we can get better constants for a \textit{large} number of $n$. To do so, we begin by defining
$$B_k := \binom{k}{\lfloor k/2\rfloor} \,=\, \left(\sqrt{\frac{2}{\pi}}+o(1)\right)\frac{2^k}{\sqrt k}$$
and $U(n):=\max\{k:\ B_k\le n\}$. The asymptote above is an application of Stirling's approximation. For a sharper two sided bound, one can use the so called Robbins' estimate \cite{robin}. Also, let $c_\ast:=\frac 12 \log_2\left(\frac\pi2\right)$.

\begin{thm}\label{density}
    For any $c_\ast-1<C<c_\ast$, the set
    $$\mathcal N_C \,:=\, \left\{n\in\mathbb N:\, U(n) \,\le\, \log_2 n + \frac 12 \log_2\log_2 n + C\right\}$$
    has asymptotic density
    $$\mathbf d\left(\mathcal N_C\right) \,=\, 2-2^{c_\ast -C}$$
    which is positive. In particular, for $C=0$, $\mathbf d\left(\mathcal N_0\right)=2-\sqrt{\frac\pi2}\approx 0.747$.
\end{thm}
\begin{proof}
    Fix a large $k$ and consider some $n\in \mathcal I_k := [B_k,B_{k+1})$, and hence $U(n)=k$. Define
    $$\rho_k \,:=\, \frac{B_{k+1}}{B_k}$$
    and write $n=B_k r$ with $r\in[1,\rho_k)$. Using $\rho_k=2(1+o(1))$, we may treat $r\in[1,2)$ asymptotically. Now, using Robbins' estimate \cite{robin}, we have
    $$\log_2 n \,=\, \log_2 B_k + \log_2 r \,=\, k - \frac12 \log_2 k - c_* + \log_2 r + o(1)$$
    implying
    \begin{align*}
        \log_2 n + \frac 12\log_2\log_2 n + C \, &=\, \left(k - \frac 12\log_2 k - c_* + \log_2 r\right) + \frac 12\log_2 k + C + o(1)\\
        &=\, k + \log_2 r + (C-c_*) + o(1)
    \end{align*}
    since $\log_2\log_2 n = \log_2 k + o(1)$ uniformly for $n\in \mathcal I_k$.

    So, the statement
    $$U(n) \,=\, k \,\le\, \log_2 n + \frac12 \log_2\log_2 n + C$$
    is equivalent to $0 \,\le\, \log_2 r + (C-c_*) + o(1) \iff r \ge 2^{c_*-C}(1+o(1))$.

    This selects the subinterval
    $$\left[B_k\cdot 2^{c_*-C}(1+o(1)),\, B_{k+1}\right)$$
    inside $\mathcal I_k$ whose relative size tends to
    $$\frac{\rho_k - 2^{c_*-C}}{\rho_k-1} \longrightarrow 2 - 2^{c_*-C}$$
    which is also the limit of the asymptotic density since the lengths $|\mathcal I_k|$ grow geometrically.
\end{proof}

\section{The Greedy Algorithm}
This section is concerned with the greedy algorithm of generating dissociated sets. It starts with two given integers $a>0$ and $b>a$. It sets $\gamma_1=a$, $\gamma_2=b$ and in the $r$-th step it chooses $\gamma_r$ to be the smallest integer greater than $\gamma_{r-1}$ such that the sequence $\{\gamma_k\}_{k=1}^r$ is dissociated.

\begin{thm}\label{greed:dissinc}
    For any given $a,b$, let $\Gamma:=\Gamma_{a,b}=\{\gamma_1=a, \gamma_2=b, \gamma_3, \dots \}$ be the dissociated sequence that the greedy algorithm produces. Then, there is an $n_0=n_0(a,b)$ such that
    $$\gamma_n = 2\cdot \gamma_{n-1}$$
    for all $n\ge n_0$.
\end{thm}

Before going to the proof, we need some definitions and a few preliminary lemmas. Define
For a finite $A\subset \mathbb N$, we define
$$\Sigma(A) \,:=\, \left\{\sum_{x\in T}x \,:\, T\subseteq A\right\}\subset\mathbb N\cup\{0\}$$
to be the family of subset sums. We also define
$$\Delta(A) \,:=\, \Sigma(A)-\Sigma(A) \,=\, \{u-v:\, u,v\in\Sigma(A)\}$$
to be the difference set of subset sums. Note that if $d=u-v\in\Delta(A)$, then $-d=v-u\in\Delta(A)$ and hence $\Delta(A)$ is symmetric. Also, for given integers $a$, $b>a$, we let $\Gamma_{a,b}=\{\gamma_1=a, \gamma_2=b, \gamma_3, \dots \}$ denote the dissociated sequence that the greedy algorithm produces. We further denote $\Gamma_{a,b}(\kappa):=\Gamma_{a,b}\cap [1,\kappa]$ for some integer $\kappa$.

We also define
$$B_n:=\left\{\sum_{i\in I}\gamma_i:\ I\subseteq[n]\right\} \,=\, \Sigma\Big(\Gamma(\gamma_n)\Big) \subset\mathbb Z$$
and
$$D_n:=\left\{\sum_{i=1}^n \varepsilon_i\gamma_i:\ \varepsilon_i\in\{-1,0,1\}\right\} \,=\, \Delta\Big(\Gamma(\gamma_n)\Big)\subset\mathbb Z$$
to be the \textit{subset-sum set} and the \textit{signed-sum set} respectively. For a finite set $A$, we denote
$$\mathscr S(A) \,:=\, \sum_{a\in A} a$$
and $\mathscr S_r := \mathscr S\Big(\Gamma(\gamma_r)\Big)$. We also let $\mathscr T_r := \mathscr S_r+1$.

\begin{lem}\label{iff}
    Fix $n\ge 1$ and $x\in\mathbb N \setminus \Gamma(\gamma_n)$. Then $\Gamma(\gamma_n)\cup\{x\}$ is dissociated if and only if $x\notin D_n$.
\end{lem}
\begin{proof}
    If $x\in D_n$, then
    $$x-\sum_{i=1}^n \varepsilon_i\gamma_i \,=\, 0$$
    for some $\varepsilon_i\in\{-1,0,1\}$ and hence $\Gamma(\gamma_n)\cup\{x\}$ is not dissociated.

    Conversely, if $\Gamma(\gamma_n)\cup\{x\}$ is not dissociated, then there exist
    $$\sum_{i=1}^n \varepsilon_i\gamma_i + \varepsilon_x x = 0$$
    for some $\varepsilon_i\in\{-1,0,1\}$ for $i\le n$ and $\varepsilon_x\in\{-1,0,1\}$, not all zero. If $\varepsilon_x=0$, this would be a forbidden nontrivial relation inside $\Gamma(\gamma_n)$, contradicting that it is dissociated. Hence $\varepsilon_x=\pm 1$ and hence
    $$x \,=\, -\varepsilon_x\sum_{i=1}^n \varepsilon_i\gamma_i \,\in\, D_n$$
    hence completing the proof.
\end{proof}

\begin{lem}\label{S_n}
    For every $n\ge 1$, $\mathscr T_n\notin D_n$. In other words, $\gamma_{n+1}\le \mathscr T_n$ for all $n\ge 2$.
\end{lem}
\begin{proof}
    This follows from the observation
    $$\left|\sum_{i=1}^n\varepsilon_i\gamma_i\right|\le \sum_{i=1}^n |\varepsilon_i|\gamma_i \le \mathscr S_n$$
    implying $D_n\subseteq [-\mathscr S_n,\mathscr S_n]$.
\end{proof}

\begin{lem}\label{T_n>2^n}
    For all $n\ge 1$, we have $\mathscr T_n\ge 2^n$.
\end{lem}
\begin{proof}
    Follows from the fact that $\Gamma(\gamma_n)$ is dissociated implying $|B_n|=2^n$ and $\Gamma$ increasing.
\end{proof}

Equipped with these, we present the promised proof.
\begin{proof}[Proof of Theorem \ref{greed:dissinc}]
    We proceed in three major steps. First, for $n\ge 1$, we define
    $$r_n \,:=\, \frac{\mathscr T_n}{2^n} \,\ge \, 1$$
    by Lemma \ref{T_n>2^n}. Also, $\{r_n\}$ is non-increasing by Lemma \ref{S_n}.

    Notice that there exists $N_\ast$ such that
    $$\gamma_{n+1} \,>\, \frac 12 \mathscr T_n$$
    for all $n\ge N_\ast$. Indeed, if not, then for infinitely many $n$, we have
    $$\mathscr T_{n+1}=\mathscr S_n+\gamma_{n+1}+1  \le \mathscr T_n+\frac12\mathscr T_n=\frac32\mathscr T_n$$
    and hence
    $$r_{n+1}=\frac{\mathscr T_{n+1}}{2^{n+1}} \le \frac{(3/2)\mathscr T_n}{2^{n+1}}=\frac34 r_n$$
    which implies $\{r_n\}$ goes to $0$ over a subsequence, a contradiction.

    The next step is to show that the greedy constraint forces large forbidden intervals. We will prove that $\{b+1,b+2,\dots,\gamma_n-1\}\subseteq D_{n-1}$ for all $n\ge 3$. We proceed by Induction. The base case is obvious since $\gamma_3$ is the smallest integer $>b=\gamma_2$ not in $D_2$.

    Now assume $\{b+1,\dots,\gamma_n-1\}\subseteq D_{n-1}$. Since $D_{n-1}\subseteq D_n$ (allowing an extra coordinate $\varepsilon_n=0$), we have $\{b+1,\dots,\gamma_n-1\}\subseteq D_n$. Also $\gamma_n\in D_n$ (take $\varepsilon_n=1$ and all other $\varepsilon_i=0$). Finally, by the greedy constraint, every integer $x$ with $\gamma_n<x<\gamma_{n+1}$ must lie in $D_n$. This completes the induction.

    Finally, fix $n\ge 3$, $1\le u\le \gamma_n-b-1$ and $t:=\gamma_n-u$. Then, $b+1\le t\le \gamma_n-1$ and hence $t,-t\in D_{n-1}$ implying $u=\gamma_n + (-t)\in \gamma_n + D_{n-1}\subseteq D_n$ (using $\varepsilon_n=1$ plus a representation of $-t$ from the first $n-1$ terms).

    This brings us to the third and final step where we prove $\gamma_{n+1}=2\gamma_n$ for all large $n$. Let's choose $n$ large enough so that $D_{n-1}$ contains all integers $1,2,\dots,b$ and $2\gamma_n> \mathscr S_{n-1}$. Now,
    $$\big[\gamma_{n-1}-b-1\big]\subseteq D_{n-1}$$
    so that, if $\gamma_{n-1}\ge 2b+1$, then $\gamma_{n-1}-b-1\ge b$ implying $[b]\subseteq D_{n-1}$. Since $\Gamma$ is strictly increasing and unbounded, there exists $N_b$ such that $\gamma_n\ge 2b+1$ for all $n\ge N_b$. Now, set $n_0:=\max\{N_b+1,\ N_\ast+1\}$. We claim that $\gamma_{n+1}=2\gamma_n$ for every $n\ge n_0$. So, from now on, fix an $n\ge n_0$.

    Since $n-1\ge N_b$, we have $[b]\subseteq D_{n-1}$ implying $\{b+1,\dots,\gamma_n-1\}\subseteq D_{n-1}$. Hence $[\gamma_n-1] \subseteq D_{n-1}$. Therefore, for each $t\in [\gamma_n-1]$, we have $\gamma_n+t\in D_n$ (taking $\varepsilon_n=1$ and representing $t\in D_{n-1}$). So, $\big\{\gamma_n+1,\gamma_n+2,\dots,2\gamma_n-1\big\}\subseteq D_n$.

    On the other hand, we will now show that $2\gamma_n\notin D_n$. If possible, let $2\gamma_n\in D_n$. Then
    $$2\gamma_n=\sum_{i=1}^n \varepsilon_i\gamma_i$$
    for $\varepsilon_i\in\{-1,0,1\}$. Let $\varepsilon:=\varepsilon_n\in\{-1,0,1\}$ and set
    $$m:=\sum_{i=1}^{n-1}\varepsilon_i\gamma_i\in D_{n-1}$$
    implying
    $$2\gamma_n=\varepsilon\gamma_n+m$$
    and hence $m=(2-\varepsilon)\gamma_n$.

    If $\varepsilon=1$, then $m=\gamma_n$. But $\gamma_n\notin D_{n-1}$ because otherwise (by Lemma \ref{iff}) adding $\gamma_n$ to $\Gamma(\gamma_{n-1})$ would not preserve the dissociated property, hence a contradiction.

    If $\varepsilon=0$, then $m=2\gamma_n$. But any $m\in D_{n-1}$ satisfies $|m|\le \mathscr S_{n-1}$ (using Lemma \ref{S_n}), while $2\gamma_n>\mathscr S_{n-1}$) since $n\ge n_0\ge N_\star+1$, a contradiction.

    If $\varepsilon=-1$, then $m=3\gamma_n$ which is again impossible because $|m|\le \mathscr S_{n-1}<2\gamma_n<3\gamma_n$.

    This completes the proof.
\end{proof}

The next natural question is to ask for an explicit $n_0(a,b)$. A clean upper bound follows essentially from the ideas discussed in the previous proof. So, we express this as a corollary.
\begin{cor}
We have
$$n_0(a,b)\ \le\ L+2K+5$$
where $L:=\left\lceil \log_2(2b+2)\right\rceil$ and $K:=\left\lceil \log_2 L\right\rceil$.
\end{cor}
\begin{proof}
Since $m\le r=L+K+2$, we have
$$n_0(a,b)\le m+1+\left\lceil\log_2 m\right\rceil$$
and hence, it is enough to bound $\lceil\log_2 m\rceil$ in terms of $K$. Since $m\le r\le L+K+2\le 4L$ (using $K\le L$ and $L\ge 2$), we have
$$\log_2 m \,\le\, \log_2(4L) \,=\, 2+\log_2 L \le 2+K$$
and thus
$$\left\lceil\log_2 m\right\rceil \le K+2$$
so that substituting $m\le L+K+2$ gives
$$n_0(a,b)\le (L+K+2)+1+(K+2)=L+2K+5$$
as claimed.
\end{proof}

\section{Generalization I}
The goal of this section and the next is to consider two different generalizations of the idea of a dissociated set.
\begin{defi}
    A set $\mathcal S =\{a_1<a_2<\dots\} \subset \mathbb N$ is said to be a $\mathcal D[g]$ set if
    $$r_{\mathcal S \cap [1,a_m]}(t) \,:=\, \left|\left\{I\subseteq[m]:\, \sum_{i\in I} a_i=t\right\}\right| \,\le\, g$$
    for all $t\in \mathbb Z$ and $m\in \mathbb N$.
\end{defi}

Under this terminology, a $\mathcal D[1]$ set is precisely the dissociated set that we have discussed before. We can upper bound $\mathcal D[g]$ sets using similar ideas as in \cite{xu}.

\begin{thm}\label{thm:dens_g}
    Let $A={a_1<\cdots<a_m}\subset\mathbb N$ be a $\mathcal D[g]$ set. Then
    $$a_m \,\ge\, \frac1g\binom{m}{\lfloor \frac m2\rfloor}-1$$
    and hence $\binom{m}{\lfloor m/2\rfloor}\le g(a_m+1)$.
\end{thm}
\begin{proof}
    Define the Boolean cube
    $$V:=\left\{-\frac 12,\frac 12\right\}^m$$
    and for $\varepsilon=(\varepsilon_1,\dots,\varepsilon_m)\in V$, define
    $$s(\varepsilon) \,:=\, \sum_{i\in I(\varepsilon)} a_i$$
    where $I(\varepsilon) \,:=\, \left\{i\in[m]:\varepsilon_i=\frac 12\right\}$. Now, take
    $$S:=\sum_{i=1}^m a_i$$
    and define the linear form
    \begin{align*}
        L(\varepsilon) \,:&=\, \sum_{i=1}^m a_i\varepsilon_i\\
        &=\frac12\sum_{i\in I(\varepsilon)}a_i-\frac12\sum_{i\notin I(\varepsilon)}a_i\\
        &=\sum_{i\in I(\varepsilon)}a_i-\frac12\sum_{i=1}^m a_i\\
        &=s(\varepsilon)-\frac S2
    \end{align*}
    and hence $L(\varepsilon)=L(\varepsilon') \Longleftrightarrow s(\varepsilon)=s(\varepsilon')$. So, each value of $L(\varepsilon)$ is attained by at most $g$ points $\varepsilon\in V$.

    Now, define
    \begin{align*}
        &V_-:= \,\{\varepsilon\in V:\ L(\varepsilon)<0\}\\
        &V_0:= \,\{\varepsilon\in V:\ L(\varepsilon)=0\}\\
        &V_+:= \,\{\varepsilon\in V:\ L(\varepsilon)>0\}
    \end{align*}
    and note that $\varepsilon\mapsto-\varepsilon$ is a bijection $V_-\leftrightarrow V_+$ implying $|V_-|=|V_+|$. Choose any subset $V_0'\subseteq V_0$ with
    $$|V_0'|=2^{m-1}-|V_-|=\left\lceil\frac{|V_0|}{2}\right\rceil$$
    and define $F:=V_-\cup V_0'$. Then, $|F|=2^{m-1}$.

    Now define the boundary
    $$\partial F \,:=\, \big\{\eta\in V\setminus F:\ \exists i\in[m]\ \text{such that }\eta \text{ differs from some }\xi\in F\text{ only in coordinate }i\big\}$$
    so that $\eta$ is outside but within Hamming distance $1$ of $F$. We have
    $$|\partial F|\ \ge\ \binom{m}{\lfloor \frac m2\rfloor}$$
    using Harper's vertex-isoperimetric inequality \cite{har}.

    Fix $\eta\in\partial F$ and fix its corresponding $\xi\in F$. Therefore, we have
    $$L(\eta)-L(\xi) \,=\, \sum_{j=1}^m a_j(\eta_j-\xi_j) \,=\, a_i(\eta_i-\xi_i) \,\in\, \{+a_i,-a_i\}$$
    since $\eta_i-\xi_i\in \{1,-1\}$. Now $L(\xi)\le 0$ and $\eta\notin F$ implying $\eta\notin V_-\cup V_0'$. So, $L(\eta)\ge 0$ with equality only if $\eta\in V_0\setminus V_0'$. So every $\eta\in\partial F$ satisfies
    $$0\le L(\eta)\le a_i\le a_m$$
    and all but those $\eta\in V_0\setminus V_0'$ satisfy $0< L(\eta)\le a_m$.

    Also, we have
    $$\left|\partial^+F\right|:= \left\lvert \{\eta\in\partial F:\ L(\eta)>0\}\right\rvert\ \ge\ |\partial F|-g$$
    since $|V_0|=r_A(S/2)\cdot \mathbf 1_{\{S \text{ is even}\}}\le g$.

    Now all values $L(\eta)$ for $\eta\in\partial^+F$ lie in the interval $(0,a_m]$, and any two distinct values differ by an integer because
    $$L(\eta)-L(\eta') \,=\, L(\eta-\eta') \,=\, \sum_{i=1}^m a_i(\eta_i-\eta'_i)$$
    and $\eta_i-\eta'_i\in\{-1,0,1\}$, so the difference is an integer. Hence the set $\{L(\eta):\eta\in\partial^+F\}\subset(0,a_m]$ has at most $a_m$ distinct values (integer gaps $\ge 1$ inside an interval of length $a_m$). With multiplicity $\le g$ per value, we have
    $$|\partial^+F| \le g\cdot a_m$$
    implying
    $$a_m \ge \frac{|\partial F|}{g}-1 \,\ge\, \frac1g\binom{m}{\lfloor \frac m2\rfloor}-1$$
    hence completing the proof.
\end{proof}

\begin{cor}\label{dens-g}
    Let $\mathcal S$ be a $\mathcal D[g]$ set. Then
    $$\left\lvert\mathcal S(n)\right\rvert \,\le\, \log_2 n \,+\, \frac 12\log_2\log_2 n \,+\, \log_2 g \,-\, \log_2\left(\sqrt{\frac 2\pi}\right) \,+\, o(1)$$
    where $\mathcal S(n):=\mathcal S\cap [1,n]$.
\end{cor}

\textit{Remark}: It should be briefly noted that there are some cute consequences of Corollary \ref{dens-g}. For instance, it immediately implies that for a given $k$ and $\ell$, there are infinitely many $N$ that can be written as the sum of distinct $k$-th powers in at least $\ell$ different ways. Similarly, there are infinitely many $N\in \mathbb N$ that can be written as the sum of primes (on your favourite arithmetic progression) in $\ell$ different ways.

On the other hand, the greedy algorithm for $\mathcal D[g]$ sets also exhibit the same doubling property as discussed, with the proof being very similar to that of Theorem \ref{greed:dissinc}. So, instead of giving the full proof again, we will only provide an outline.

\begin{thm}\label{greedy-Dg}
    For any given $a,b$, let $\gamma_1=a$, $\gamma_2=b$ and in let $\gamma_r$ to be the smallest integer greater than $\gamma_{r-1}$ such that the sequence $\{\gamma_k\}_{k=1}^r$ is $\mathcal D[g]$. Let $\Gamma:=\Gamma_{a,b}^{(g)}=\{\gamma_1=a, \gamma_2=b, \gamma_3, \dots \}$. Then, there is an $n_0=n_0(a,b,g)$ such that
    $$\gamma_n = 2\cdot \gamma_{n-1}$$
    for all $n\ge n_0$.
\end{thm}
\begin{proof}
    As advertised, the proof is formally identical to the dissociated ($g=1$) proof once one replaces the signed–sum set $D_n$ by the forbidden–shift condition expressed via the representation function
    $$r_n(s) \,:=\, \left|\left\{I\subseteq [n] \,:\, \sum_{i\in I}\gamma_i=s,\right\}\right| \,=\, \big[z^s\big] \prod_{i=1}^n\big(1+z^{\gamma_i}\big)$$
    for $s\in \mathbb Z$.

    Adding $x$ transforms $r_n$ by the standard subset-sum recurrence
    $$r_{n+1}(s)=r_n(s)+r_n(s-x)$$
    and hence $\Gamma\big(\gamma_n\big)\cup\{x\}$ is $\mathcal D[g]$ iff $r_n(s)+r_n(s-x)\le g$ for all $s$ iff $x$ is not in the forbidden set
    $$F_n \,:=\, \left\{x\in\mathbb N:\exists s\ \text{with }r_n(s)+r_n(s-x)\ge g+1\right\}$$
    which will act as the analogue of $D_n$ from the previous proof.

    Finally, borrowing notations from the previous proof, we have
    $$2^n \,=\, \sum_s r_n(s) \,\le\, g(\mathscr S_n+1) \,=\, g\mathscr T_n$$
    so the potential becomes
    $$r_n^{(g)} \,:=\, \frac{g\mathscr T_n}{2^n}$$
    and the rest of the proof runs through.
\end{proof}

\section{Generalization II}
\begin{defi}
    A set $\mathcal S\subset \mathbb N$ is called a $\mathcal D_k$ set if any equality of the form
    $$\sum_{s\in \mathcal S} \varepsilon_s \cdot s =0$$
    where $\varepsilon_s \in \{-k,\dots ,-1,0,+1, \dots ,k\}$ implies that all the $\varepsilon_s$'s are $0$.
\end{defi}

Dissociated sets in this language are $\mathcal D_1$ sets and they have been called \textit{detecting sets} in literature \cite{diss4, detect}. In this case again, we can get an upper bound by generalizing the argument in \cite{xu}. This is what we will now exhibit.

\begin{thm}\label{thm:dens-k}
    Let $\mathcal S=\{a_1<\cdots<a_m\}\subset\mathbb N$ be a $\mathcal D_k$ set. Then, we have
    $$a_m \,\ge\, (1+o(1))\, \frac{(k+1)^m}{2}\, \sqrt{\frac{6}{m\pi k\, (k+2)}}$$
    and hence, for an infinite $\mathcal D_k$ set $\mathcal S$, we have
    $$\left\lvert\mathcal S(n)\right\rvert \,\le\, \log_{k+1} n \,+\, \frac 12 \log_{k+1}\log_2 n \,+\, \frac 12 \log_{k+1} \left(\frac{2\pi}{3}\,k(k+2)\right) \,+\, o(1)$$
    where $\mathcal S(n) \,:=\, \mathcal S\cap [1,n]$.
\end{thm}
\begin{proof}
    As before, we again define
    $$V:=[0,k]^m=\{0,1,\dots,k\}^m$$
    viewed as the vertex set of the grid graph
    $$G= \mathop{\mathlarger{\mathlarger{\mathlarger{\times}}}}_{i\in [m]} P_{k+1}$$
    expressed as the Cartesian product of path graphs - two vectors are adjacent iff they differ by $\pm 1$ in exactly one coordinate. For $x=(x_1,\dots,x_m)\in V$, define
    $$s(x):=\sum_{i=1}^m a_i x_i\in\mathbb Z,\quad S:=\sum_{i=1}^m a_i, \quad L(x):=s(x)-\frac{k}{2}S$$
    and note $L(k\cdot \mathbf 1-x)=-L(x)$ where $\mathbf 1=(1,\dots,1)$. Let
    $$c_{m,k}(t) \,:=\, \big[x^t\big]\,\Big(1+x+\cdots+x^k\Big)^m$$
    be the number of vectors $x\in V$ with $|x|:=x_1+\cdots+x_m=t$. The sequence $t\mapsto c_{m,k}(t)$ is symmetric and unimodal, and it attains its maximum at $t=\lfloor \frac{km}{2}\rfloor$.

    Now, note that $\mathcal S$ is a $\mathcal D_k$ set, then 
    $$s(x)=s(y) \implies 0=s(x)-s(y)=\sum_{i=1}^m (x_i-y_i)a_i$$
    with $x_i-y_i\in\{-k,\dots,k\}$ and hence, the map $s:V\to\mathbb Z$, $x\mapsto \sum a_ix_i$ is injective. In particular, $L$ is also injective.

    Again as before, we set
    \begin{align*}
        &V_- \,:=\, \{x\in V:L(x)<0\}\\
        &V_0 \,:=\, \{x\in V:L(x)=0\}\\
        &V_+ \,:=\, \{x\in V:L(x)>0\}
    \end{align*}
    and use $L(k\mathbf\cdot 1-x)=-L(x)$ to get $|V_-|=|V_+|$.

    Again, choose $V_0'\subseteq V_0$ so that
    $$\big|F\big| \,:=\, \big|V_-\big| \,+\, \big|V_0'\big| \,=\, \left\lfloor\frac{|V|}{2}\right\rfloor$$
    and hence, $F\subseteq V$ has size $\left\lfloor \frac{(k+1)^m}2 \right\rfloor$. Finally, also define the vertex boundary $\partial F$ to be the set of all $y\in V\setminus F$ such that there is $x\in F$ adjacent to $y$ in the grid.

    From Theorem 11 in \cite{leader} (also see \cite{wang}), we have
    $$|N(F)\setminus F| \,=\, |N(F)|-|F| \,\ge\, |N(C)|-|C|$$
    implying $|\partial F| \,\ge\, |\partial C|$.

    On the other hand, let $t_0:= \left\lfloor \frac{km}2\right\rfloor$ and $L_{t}:=\{x\in V:|x|=t\}$) so that $|L_t|=c_{m,k}(t)$. Because the level sizes are symmetric and unimodal with maximum at $t_0$ (see \cite{li} for example), we have
    $$\sum_{t< t_0} \big|L_t\big| \,\le\, \frac{\big|V\big|-\big|L_{t_0}\big|}{2}$$
    so that when you take any initial segment $C$ of size $\left\lfloor \frac{|V|}2\right\rfloor$, you must include at most half of the middle layer $L_{t_0}$, i.e., $|C\cap L_{t_0}|\le \frac{|L_{t_0}|}{2}$.

    But, every vertex in $L_{t_0}\setminus C$ has a neighbor in $L_{t_0-1}\subseteq C$ and hence
    $$L_{t_0}\setminus C \,\subseteq\, \partial C$$
    implying
    $$\big|\partial F\big| \,\ge\, \big|\partial C\big| \,\ge\, \big|L_{t_0}\setminus C\big| \,\ge\, \big|L_{t_0}\big|-\frac{|L_{t_0}|}{2} \,=\, \frac12\,\big|L_{t_0}\big| \,=\, \frac12\, M_{m,k}$$
    where $M_{m,k}$ is the maximum of $c_{m,k}(t)$ over $t$.

    Now, consider any $y\in \partial F$ with $L(y)>0$. By definition of boundary, there exists an adjacent $x\in F$ differing from $y$ in exactly one coordinate $i$ by $\pm 1$. Since $x\in F$ implies $L(x)\le 0$ and $L(y)>0$, necessarily $L(y)-L(x)=+a_i$. Hence
    $$0 \,<\, L(y) \,=\, L(x)+a_i \,\le\, a_i \,\le\, a_m$$
    and hence every such $y$ yields a value $L(y)\in(0,a_m]$.

    Moreover, for any $y,y'\in V$, we have
    $$L(y)-L(y')=s(y)-s(y')\in\mathbb Z$$
    and hence all $L$-values lie in a single coset of $\mathbb Z$ and are spaced by integers. Since $L$ is injective, distinct $y$’s give distinct $L(y)$’s. Therefore the set
    $$\big\{L(y):y\in\partial F,\ L(y)>0\big\}\subset(0,a_m]$$
    has size at most $a_m$. Finally, because $|V_0|\le 1$, at most one boundary vertex can have $L(y)=0$. Hence
    $$a_m \,\ge\, \big|\partial F\big|-1 \,\ge\, \frac12\, M_{m,k}-1$$
    since $|\partial F|\,\le\, a_m+1$.

    By unimodality of $M_{m,k}$ and Central Limit Theorem, we have (see \cite{li}and \cite{xu})
    $$M_{k,m}\sim (k+1)^m\sqrt{\frac{6}{m\pi k\, (k+2)}}$$
    and equivalently, the $(k+1)$-multinomial maximum $\sup_t \binom{m,k+1}{t}$ satisfies the same estimate. This completes the proof.
\end{proof}

Of course, the next question is to ask about the greedy algorithm for $\mathcal D_k$ sets. This is what we will now provide. The proof of the following theorem requires three lemmas that we will provide after giving the proof.

\begin{thm}\label{greedy-Dk}
    Fix an integer $k\ge 1$. Let $a\in\mathbb N$, let $\gamma_1:=a$ and let $\gamma_r$ be the smallest integer greater than $\gamma_{r-1}$ such that $\{\gamma_i\}_{i=1}^r$ is a $\mathcal D_k$ set. Then there exists $n_0=n_0(a,k)$ so that
    $$\gamma_n=(k+1)\cdot \gamma_{n-1}$$
    for all $n\ge n_0$.
\end{thm}
\begin{proof}
    Again write $\Gamma(\gamma_n):=\{\gamma_1,\dots,\gamma_n\}$ and define
    $$\mathscr S_n:=\sum_{i=1}^n \gamma_i, \qquad \mathscr T_n:=\ k\,\mathscr S_n+1$$
    as before. Also write
    $$B_n:=\left\{\sum_{i=1}^n x_i\gamma_i :\ x_i\in\{0,1,\dots,k\}\right\}\subset \big[0,k \mathscr S_n\big] = \big[0, \mathscr T_n-1\big]$$
    and $D_n:=B_n-B_n\subset \big[-(\mathscr T_n-1),\, \mathscr T_n-1 \big]$. If $\Gamma(\gamma_n)$ is a $\mathcal D_k$ set, then the map
    $$\{0,1,\dots,k\}^n\to\mathbb Z,\qquad (x_i)\mapsto \sum_{i=1}^n x_i\gamma_i$$
    is injective, implying $|B_n|=(k+1)^n\le \mathscr T_n$.
    
    Also, define the \textit{slack}
    $$\delta_n:= \mathscr T_n-\gamma_{n+1}\in\mathbb Z_{\ge 0}$$
    and the normalized potential
    $$r_n:=\frac{\mathscr T_n}{(k+1)^n}$$
    which satisfies $r_{n+1}\le r_n$ using
    $$\mathscr T_{n+1}= k \mathscr S_{n+1}+1= k(\mathscr S_n+\gamma_{n+1})+1= \mathscr T_n+k\gamma_{n+1}$$
    and Lemma \ref{2.2} proved below. So $r_n$ is monotone non-increasing and bounded below by $1$, implying $r_\infty:=\lim r_n\in[1,\infty)$.

    Now, divide the identity $\mathscr T_{n+1}= \mathscr T_n+k\gamma_{n+1}$ by $(k+1)^{n+1}$ to get
    $$r_{n+1}=\frac{r_n}{k+1}+\frac{k}{k+1}\cdot \frac{\gamma_{n+1}}{(k+1)^n}$$
    implying
    $$\frac{\gamma_{n+1}}{(k+1)^n}=\frac{(k+1)r_{n+1}-r_n}{k}$$
    implying
    $$\frac{\gamma_{n+1}}{(k+1)^n}\longrightarrow \frac{(k+1)r_\infty-r_\infty}{k}=r_\infty$$
    using $r_n,r_{n+1}\to r_\infty$. Therefore, we have
    $$\frac{\gamma_{n+1}}{\gamma_n} =(k+1)\cdot \frac{\gamma_{n+1}(k+1)^{n-1}}{\gamma_n \ (k+1)^{n}} \longrightarrow (k+1)\cdot \frac{r_\infty}{r_\infty}=k+1$$
    implying
    $$\frac{\delta_n}{(k+1)^n} =\frac{\mathscr T_n-\gamma_{n+1}}{(k+1)^n} \longrightarrow r_\infty-r_\infty=0$$
    since $\gamma_{n+1}$ is the smallest admissible integer $>\gamma_n$, it must satisfy $\gamma_{n+1}\le x=\mathscr T_n$
    since $\mathscr T_n\to r_\infty (k+1)^n$ and $\gamma_{n+1}\to r_\infty (k+1)^n$.

    First, we claim that there exists $N_*$ such that for all $n\ge N_*$, $2\gamma_{n+1}> \mathscr T_n$. Indeed, if possible, let $2\gamma_{n+1}\le \mathscr T_n$ for infinitely many $n$. For such $n$, we have
    $$\mathscr T_{n+1}= \mathscr T_n+k\gamma_{n+1}\le \mathscr T_n+\frac{k}{2} \mathscr T_n=\frac{k+2}{2}\, \mathscr T_n$$
    implying $2(k+1)\, r_{n+1}\le (k+2)\, r_n$. Let $c:=\frac{k+2}{2(k+1)}<1$. If the inequality $r_{n+1}\le c r_n$)\ holds infinitely often, then after $m$ such occurrences we would have $r_n\le c^m\, r_{n_0}$, which goes to $0$ as $m\to\infty$. This is contradiction, hence completing the proof.

    Fix $n\ge N_*$ and put $M:=\gamma_{n+1}$. Define
    $$\mathscr H_n:=\left\lceil \frac{\mathscr T_n}{2}\right\rceil,\quad \mathscr L_n:=\max\bigl\{\gamma_n+1,\ \mathscr H_n\bigr\}$$
    and note that $\mathscr L_n\le M$ for all $n\ge N_*$.

    Now, we claim that for all $n\ge N_*$, we have $\big[\mathscr L_n,\ M-1\big]\cap\mathbb Z \subset D_n$. Indeed, take any integer $y$ with $\mathscr L_n\le y\le M-1$. Then $y>\gamma_n$ and $y<M=\gamma_{n+1}$. By greediness, $\gamma_{n+1}$ is the smallest admissible integer $>\gamma_n$, so every integer strictly between $\gamma_n$ and $\gamma_{n+1}$ is not admissible at stage $n$. Hence $y$ is forbidden at stage $n$ implying there exists $d\in\{1,\dots,k\}$ with $dy\in D_n$. But $y\ge \mathscr H_n$ implying $2y> \mathscr T_n-1$. Therefore for every $d\ge 2$, $dy\ge 2y> \mathscr T_n-1$, implying $D_n\subset\big[-(\mathscr T_n-1),\, \mathscr T_n-1\big]$, so $dy\notin D_n$ for all $d\ge 2$. Hence the only possible witness is $d=1$, so $y\in D_n$. This holds for every $y\in\big[\mathscr L_n,\, M-1\big]$, hence completing the proof.

    Now, we have
    $$B_{n+1}=B_n+\big\{0,M,2M,\dots,kM\big\}$$
    implying
    $$D_{n+1}=B_{n+1}-B_{n+1}=D_n+ \big\{-kM,-(k-1)M,\dots,0,\dots,(k-1)M,kM\big\}$$
    and in particular, $D_n+kM\subset D_{n+1}$ and finally, $\big[kM+\mathscr L_n,\ kM+(M-1)\big] = \big[kM+\mathscr L_n,\ (k+1)M-1\big] \subset D_{n+1}$. Also, from Lemma \ref{7.1} proved below, we have $\gamma_{n+2}\le (k+1)M$.

    Finally, assume $n\ge N_*$. If $\gamma_{n+2}<(k+1)M$, then $\gamma_{n+2}\le (k+1)M-1$. But every integer $x\in\big[kM+\mathscr L_n,\, (k+1)M-1\big]$ belongs to $D_{n+1}$, hence is forbidden at stage $n+1$. Since $\gamma_{n+2}$ is admissible, it cannot lie in that interval. Therefore, we have
    $$\frac{\gamma_{n+2}}{M}\le k + \frac{\mathscr L_n}{M}$$
    since $\gamma_{n+2}\le kM+\mathscr L_n-1$. By definition, we have $\mathscr T_n=M+\delta_n$, so that
    $$\frac{\mathscr H_n}{M}=\frac{1}{M}\left\lceil \frac{\mathscr T_n}2\right\rceil =\frac{1}{2}\cdot\frac{\mathscr T_n}{M}+o(1) =\frac12\left(1+\frac{\delta_n}{M}\right)+o(1)$$
    implying
    $$\frac{\mathscr H_n}{M}\longrightarrow \frac12$$
    since $\delta_n=o(M)$. Also,we have
    $$\frac{\gamma_n}{M}=\frac{\gamma_n}{\gamma_{n+1}}\longrightarrow \frac{1}{k+1}\le \frac12$$
    implying
    $$\frac{\mathscr L_n}{M}=\max\left(\frac{\gamma_n+1}{M},\ \frac{\mathscr H_n}{M}\right)\longrightarrow \max\left(\frac{1}{k+1},\ \frac12\right)=\frac12$$
    so that there exists $N_2$ such that for all $n\ge N_2$, we have $8\mathscr L_n<5M$.

    Let $n\ge n_0:=\max \big(N_*,N_1,N_2\big)$. If $\gamma_{n+2}<(k+1)M$, then $$\frac{\gamma_{n+2}}{M}\le\ k+\frac{\mathscr L_n}{M}<\ k+\frac58$$
    which is a contradiction. We conclude that for all $n\ge n_0$, we have
    $$\gamma_{n+2}=(k+1)\gamma_{n+1}$$
    hence completing the proof.
\end{proof}

\begin{lem}\label{2.1}
    Assume $\Gamma(\gamma_n)$ is $\mathcal D_k$. For an integer $x>\gamma_n$, the set $\Gamma(\gamma_n)\cup\{x\}$ is $\mathcal D_k$ if and only if $dx\notin D_n$ for all $d\in\{1,2,\dots,k\}$.
\end{lem}
\begin{proof}
    Consider
    $$B_{n+1}=\left\{\sum_{i=1}^n x_i\gamma_i+t x :\ x_i\in\{0,\dots,k\},\ t\in\{0,\dots,k\}\right\} =\ \bigcup_{t=0}^k \big(B_n+t x\big)$$
    and note that the set $\Gamma_n\cup\{x\}$ is $\mathcal D_k$ iff the digit map on $\{0,\dots,k\}^{n+1}$ is injective, i.e., iff the union $\bigcup_{t=0}^k \big(B_n+t x\big)$ is disjoint.

    Two translates $B_n+s x$ and $B_n+t x$ intersect iff there exist $u,v\in B_n$ with $u+s x=v+t x$, implying $(t-s)x=u-v\in D_n$. Since $D_n$ is symmetric ($D_n=-D_n$), this is equivalent to requiring $dx\notin D_n$ for all $d\in\{1,\dots,k\}$.
\end{proof}

\begin{lem}\label{2.2}
    For every $n\ge 1$, the integer $\mathscr T_n = k \mathscr S_n+1$ is admissible at stage $n$. In particular, $\gamma_{n+1}\le \mathscr T_n$.
\end{lem}
\begin{proof}
    Since $D_n\subset \big[-k \mathscr S_n,\, k \mathscr S_n\big]$, if $x=\mathscr T_n= k\mathscr S_n+1$, then
    $$dx\ge x=k\mathscr S_n+1>k\mathscr S_n\ge \max D_n$$
    for every $d\in\{1,\dots,k\}$. So, $dx\notin D_n$. By Lemma \ref{2.1}, $x$ is admissible.
\end{proof}

\begin{lem}\label{7.1}
    The integer $(k+1)M$ is admissible at stage $n+1$ for every $n\ge N_*$.
\end{lem}
\begin{proof}
    We must verify
    $$d(k+1)M\notin D_{n+1}$$
    for all $d\in\{1,\dots,k\}$. We split into cases.

    \textbf{Case I} ($d=1$): Suppose $(k+1)M\in D_{n+1}$. Then, there exists $u\in D_n$ and $t\in\{-k,-k+1,\dots,k\}$ such that
    $$(k+1)M = u + tM \quad\Longrightarrow\quad u = (k+1-t)M$$
    implying $(k+1-t)M\in D_n$. Note that $k+1-t$ is an integer and since $t\in[-k,k]$, it satisfies $k+1-t\in\{1,2,\dots,2k+1\}$. If $k+1-t=1$, then $t=k$ and we get $M\in D_n$. But $M=\gamma_{n+1}$ is admissible at stage $n$, so Lemma \ref{2.1} implies $1\cdot M\notin D_n$, which is a contradiction. Also, if $k+1-t\ge 2$, then $(k+1-t)M\ge 2M$. Since $n\ge N_*$, we have $M>\frac 12 \mathscr T_n$, hence $2M>\mathscr T_n$ and so $2M> \mathscr T_n-1$. But $D_n\subset\big[-(\mathscr T_n-1),\, \mathscr T_n-1\big]$ so that no integer $\ge 2M$ lies in $D_n$, a contradiction. Hence, $(k+1)M\notin D_{n+1}$.

    \textbf{Case II} ($2\le d\le k$): We have $D_{n+1}\subset \big[-k\mathscr S_{n+1},\,k\mathscr S_{n+1}\big]$, where $\mathscr S_{n+1}= \mathscr S_n+M$. But $k \mathscr S_{n+1}=k \mathscr S_n+kM = (\mathscr T_n-1)+kM$. Since $n\ge N_*$, we have $\mathscr T_n-1<2M-1$, so that $k \mathscr S_{n+1} < (2M-1)+kM = (k+2)M-1$. On the other hand, for $d\ge 2$, we have $d(k+1)M\ge 2(k+1)M \ge k \mathscr S_{n+1}$. This implies $d(k+1)M\notin D_{n+1}$ for all $d\in\{2,\dots,k\}$.

    This completes the proof.
\end{proof}

Of course one could give similar generalizations of Theorem \ref{thm:dens-k} or Corollary \ref{dens-g} along the lines of Theorem \ref{density}. Also, the question of finding $n_0(a,b,g)$ from Theorem \ref{greedy-Dg} or that of finding $n_0(a,k)$ from Theorem \ref{greedy-Dk} can be handled similarly as done after Theorem \ref{greed:dissinc}. The details are left to the interested reader.

\section{Dissociated Sets Avoiding Geometric Progression}
Finally, we will conclude by making an observation. From Theorems \ref{greed:dissinc} and \ref{dens-g}, it might look like dissociated sets must contain lots of geometric progressions. We will show that this is probably \textit{not} true by proving it in a \textit{large} collection of dissociated sets.

To begin with, let $\mathscr D(n)$ be the collection of all dissociated subsets of $[n]$. From Theorem \ref{thm:dens}, for $M:=\lfloor 2\log_2 n\rfloor$, we have
$$\mathscr D(n) \,\le\, \sum_{j=0}^{M}\binom{n}{j} \,\ll\, (M+1)\binom{n}{M} \,\ll\, (M+1)\left(\frac{en}{M}\right)^M$$
implying $\mathscr D(n)\le \exp\Big(\mathcal O\big((\ln n)^2\big)\Big)$.

Now, fix a large $n$ and consider the following family. Set $\tilde m:=\lfloor \log_2 n\rfloor$ and $m:= \tilde m-L$ for an integer $L=\mathcal O(\log_2\log_2 n)$. For each $t\in\{0,1,\dots,m-1\}$, define the dyadic interval
$$I_t := \left(\frac{n}{2^{t+1}},\,\frac{n}{2^{t}}\right]\cap \big(2\mathbb Z+1\big)$$
and define $\mathcal F_L(n)$ to be the family of sets
$$S(\mathbf q) \,=\, \big\{2^t q_t : 0\le t\le m-1\big\}\subset [n]$$
where $q_t\in I_t$ for each $t$. Now, let
$$\mathscr G_{\mathcal F}(n):= \#\big\{S\in\mathcal F_L(n): \nexists \,a,b,c\in S : b^2=ac\big\}$$
and $\mathscr D_{\mathcal F}(n):=|\mathcal F_L(n)|$.

Now, every $\mathcal S\in \mathcal F_L(n)$ is dissociated and $|\mathcal S|\sim (\log_2 n)$. Indeed, if
$$\sum_{t=0}^{m-1}\varepsilon_t\big(2^t q_t\big)=0$$
with $\varepsilon_t\in\{-1,0,1\}$, then
$$\varepsilon_{t_0}\, 2^{t_0}\, q_{t_0}\, \equiv\,  0 \pmod{2^{t_0+1}}$$
where $t_0$ be the smallest index with $\varepsilon_{t_0}\neq 0$. But $q_{t_0}$ is odd, so
$$\varepsilon_{t_0}2^{t_0}q_{t_0}\equiv \pm 2^{t_0}\not\equiv 0\pmod{2^{t_0+1}}$$
which is a contradiction, proving $\varepsilon_t=0$ for all $t$.

On the other hand, we also have
$$\log_2 \big|\mathcal F_L(n)\big| \,=\, \sum_{t=0}^{m-1} \log_2 \big|I_t\big| \,=\, \sum_{t=0}^{m-1}\Big(\log_2 n-t+O(1)\Big) \,=\, \frac{m^2}{2} + \mathcal O\big(m\log_2 m\big)$$
implying
$$\big|\mathcal F_L(n)\big| \,=\, 2^{\frac 12(\log_2 n)^2 + o\Big((\log n)^2\Big)}$$
and hence, the family $\mathcal F_L(n)$ is \textit{quite large}. Notice that this also proves $\mathscr D(n) = \exp\Big(\Theta\big((\ln n)^2\big)\Big)$.

Now, write $x_t:=2^t q_t$ and notice that if $x_j^2=x_ix_k$, then $2j=i+k$ and $q_j^2=q_iq_k$. This implies that if we set $i,j,k$ with $i+k=2j$, then
$$\mathbb P\big(x_j^2=x_ix_k\big) \,\le\, \frac{1}{|I_k|}$$
since at most one choice of $q_k\in I_k$ can satisfy $q_kq_i=q_j^2$.

Now, for a fixed $k$, we have $\#\big\{(i,j): i<j<k,\ i+k=2j\big\} \,\le\, k/2$. Therefore, for $k\le m-1$, we have $|I_k|\gg n\,2^{-k}$. So, if $\mathcal S$ is chosen from $\mathcal F_L(n)$ uniformly randomly, then
$$\mathbb P\big(\exists\, a,b,c\in \mathcal S : b^2=ac\big) \,\ll\, \sum_{k=2}^{m-1} \left(\frac{k}{2}\right)\frac{2^k}{n} \ll \tilde m\,2^{-L}$$
and consequently
$$\frac{\mathscr G_{\mathcal F}(n)}{\mathscr D_{\mathcal F}(n)} \,\ge\, 1 - \mathcal O\Big(\tilde m\,2^{-L}\Big)$$
thus showing that most members of $\mathcal F_L(n)$ don't contain a three term geometric progression.

\textit{Remark}: We can sharpen the result to
$$\frac{\mathscr G_{\mathcal F}(n)}{\mathscr D_{\mathcal F}(n)} \,\ge\, 1 - \mathcal O \left(\frac{1}{\big(\log_2 n\big)^{A+1}}\right)$$
by choosing $L := \left\lceil (A+2)\log_2 m \right\rceil$ and doing a similar calculation.

See \cite{dutta} for a nicer explicit construction of a dissociated sequence avoiding a geometric progression.

\section{Acknowledgements}
I would like to thank Prof. Ramachandran Balasubramanian, Prof. Greg Martin, Prof. Ilya Shkredov, Prof. Sándor Kiss, Satvik Saha and Sohom Gupta for the discussions I had with them about this topic.

\bibliographystyle{plain}
\bibliography{diss}

\end{document}